\documentclass[12pt]{article}

\input diagram.tex

\usepackage{amsmath} 
\usepackage[mathscr]{eucal}
\usepackage{amssymb} 
\usepackage{theorem}
\usepackage[dvipsnames]{xcolor}

\theoremstyle{change} 
\newtheorem{theorem}{Theorem}[section] 
\newtheorem{lemma}[theorem]{Lemma} 
\newtheorem{proposition}[theorem]{Proposition}
\newtheorem{corollary}[theorem]{Corollary}
\theorembodyfont{\rmfamily} 

\newtheorem{remark}[theorem]{Remark}

\newtheorem{definition}[theorem]{Definition}
\newtheorem{notation}[theorem]{Notation}
\newtheorem{nothing}[theorem]{} 

\newtheorem{remark and notation}[theorem]{Remark and Notation}

\newenvironment{proof}{\noindent{\bf Proof}\ }{\qed\bigskip}

\renewcommand{\le}{\leqslant}
\renewcommand{\ge}{\geqslant} 

\newcommand{\abar}{\bar{a}}

\newcommand{\Aut}{\mathrm{Aut}}

\newcommand{\br}{\mathrm{br}}

\newcommand{\calF}{\mathcal{F}}

\newcommand{\calO}{\mathcal{O}}
\newcommand{\calOtilde}{\tilde{\calO}}

\newcommand{\catfont}{\mathsf}

\newcommand{\ebar}{\bar{e}}

\newcommand{\FF}{\mathbb{F}}

\newcommand{\Hom}{\mathrm{Hom}}
\newcommand{\id}{\mathrm{id}}

\newcommand{\Ind}{\mathrm{Ind}}
\newcommand{\Inf}{\mathrm{Inf}}

\newcommand{\lexp}[2]{\setbox0=\hbox{$#2$} \setbox1=\vbox to
                 \ht0{}\,\box1^{#1}\!#2}
\newcommand{\lmod}[1]{\llap{\phantom{|}}_{#1}\catfont{mod}}

\newcommand{\ltriv}[1]{\llap{\phantom{|}}_{#1}\catfont{triv}}

\newcommand{\myiso}{\buildrel\sim\over\to}

\newcommand{\qed}{\nobreak\hfill
                   \vbox{\hrule\hbox{\vrule\hbox to 5pt
                   {\vbox to 8pt{\vfil}\hfil}\vrule}\hrule}}
\newcommand{\QQ}{\mathbb{Q}}

\newcommand{\Res}{\mathrm{Res}}
\newcommand{\res}{\mathrm{res}}

\newcommand{\scrS}{\mathscr{S}}
\newcommand{\scrStilde}{\tilde{\scrS}}
\newcommand{\scrT}{\mathscr{T}}
\newcommand{\scrX}{\mathscr{X}}
\newcommand{\scrY}{\mathscr{Y}}

\newcommand{\Vtilde}{\tilde{V}}
\newcommand{\ZZ}{\mathbb{Z}}

\title{The orthogonal unit group of the trivial source ring\footnote{{\bf MR Subject Classification:}  
20C20, 20C15, 19A22. {\bf Keywords:}  Trivial source ring, $p$-permutation modules, orthogonal units, $p$-permutation equivalences, species}}

\author{\small Robert Boltje\\
  \small Department of Mathematics\\
  \small University of California\\
  \small Santa Cruz, CA 95064\\
  \small USA\\
  \small boltje@ucsc.edu
  \and
  \small Rob Carman\\
  \small Department of Mathematics\\
  \small William \& Mary\\
  \small Williamsburg, VA 23185\\
  \small USA\\
  \small wrcarman@wm.edu}
\date{October 24, 2022}

\begin{document}
\sloppy


\maketitle


\begin{abstract}
Let $G$ be a finite group, $p$ a prime, and $(K,\mathcal{O},F)$ a $p$-modular system. We prove that the trivial source ring of $\mathcal{O} G$ is isomorphic to the ring of {\em coherent} $G$-stable tuples $(\chi_P)$, where $\chi_P$ is a virtual character of $K[N_G(P)/P]$, $P$ runs through all $p$-subgroups of $G$, and the coherence condition is the equality of certain character values. We use this result to describe the group of orthogonal units of the trivial source ring as the product of the unit group of the Burnside ring of the fusion system of $G$ with the group of coherent $G$-stable tuples $(\varphi_P)$ of homomorphisms $N_G(P)/P\to F^\times$. The orthogonal unit group of the trivial source ring of $\mathcal{O} G$ is of interest, since it embeds into the group of $p$-permutation autoequivalences of $\mathcal{O} G$.
\end{abstract}


\section{Introduction}\label{sec intro}

Throughout this paper we fix a finite group $G$, a prime $p$, and a $p$-modular system $(K,\calO, F)$. In other words, $\calO$ is a complete discrete valuation ring with residue field $F$ of characteristic $p$ and field of fractions $K$ of characteristic zero.

\smallskip
The class of $p$-permutation $\calO G$-modules and its Grothendieck ring $T(\calO G)$, the trivial source ring of $\calO G$, play an important role in the modular representation theory of finite groups. For instance, several equivalence relations between block algebras (as for instance splendid Rickard equivalences, see \cite{Rickard1996}, and $p$-permutation equivalences, see \cite{BX2008} and \cite{BoltjePerepelitsky2020}) are defined in terms of $p$-permutation bimodules or their chain complexes. They give satisfactory explanations for the phenomenon of {\em isotypies}, introduced by Brou\'e, see \cite{Broue1990}. Finding new ways to describe $p$-permutation modules or the ring $T(\calO G)$ is therefore of interest. In this paper we give a  description of $T(\calO G)$ in terms of tuples of {\em coherent characters}. More precisely, we consider a ring homomorphism
\begin{equation*}
   \beta_G\colon T(\calO G)\to \Bigl(\prod_{P\in\scrS_p(G)} R(K[N_G(P)/P])\Bigr)^G\,,
\end{equation*}
where $\scrS_p(G)$ denotes the set of $p$-subgroups of $G$, $R(K[N_G(P)/P])$ denotes the ring of characters of $K[N_G(P)/P]$-modules, and the exponent $G$ stands for taking $G$-fixed points under the natural $G$-conjugation action on the product. See Section~\ref{sec Thm A} for a precise definition of the homomorphism $\beta_G$, using the Brauer construction for each $P\in\scrS_p(G)$.

\bigskip\noindent
{\bf Theorem A}\quad{\it
The ring homomorphism $\beta_G$ is injective and its image consists of those tuples $(\chi_P)\in\bigl(\prod_{P\in\scrS_p(G)} R(K[N_G(P)/P])\bigr)^G$ which satisfy the following coherence condition:
\begin{equation*}
   \text{{\rm (C)}\quad For each $P\in\scrS_p(G)$ and each $x\in N_G(P)$ one has $\chi_P(xP) = 
   \chi_{P\langle x_p\rangle}(xP\langle x_p\rangle)$.}
\end{equation*}
}

\noindent
Here, $x_p$ denotes the $p$-part of the element $x\in G$ and $\langle x_p\rangle$ the subgroup generated by $x_p$. We call $G$-fixed tuples $(\chi_P)$ of virtual characters satisfying the condition (C) {\em coherent tuples}. Theorem~A is proved in Section~\ref{sec Thm A}. It has also a block-wise version, see Corollary~\ref{cor Thm A blockwise} and a version related to limiting the set of vertices, see Corollary~\ref{cor Thm A with scrX}. These versions will be crucial when analyzing $p$-permutation equivalences between blocks in terms of the associated coherent character tuples and studying the precise relationship between isotypies and $p$-permutation equivalences, as anticipated in a future paper.

\bigskip
We will present several applications of the above theorem. A crucial tool in these applications is Theorem~B. Although it follows immediately from \cite[Theorem~7.1]{BarsottiCarman2020} (see Theorem~\ref{thm BaCa}), we think it is worth stating it in the introduction. Let $S$ denote a Sylow $p$-subgroup of $G$ and let $\calF:=\calF_S(G)$ denote the resulting fusion system on $S$. The Burnside ring $B(\calF)$ of $\calF$ is a subring of the Burnside ring $B(S)$, see \ref{def B(calF)} for a precise definition. One has a natural ring homomorphism $\lambda_G\colon B(G)\to T(\calO G)$ induced by mapping a finite $G$-set $X$ to its {\em $\calO$-linearization}, the permutation $\calO G$-module $\calO X$. Note that Green's theory of vertices and sources and Green's indecomposability theorem imply that $\lambda_S\colon B(S)\to T(\calO S)$ is a ring isomorphism, since $S$ is a $p$-group. It is an easy verification that the images of the restriction maps
\begin{equation*}
   \res^G_S\colon B(G)\to B(S)\quad \text{and}\quad \lambda_S^{-1}\circ\res^G_S\colon T(\calO G)\to T(\calO S)\myiso B(S)
\end{equation*}
are contained in $B(\calF)$. Theorem~B states that the resulting restriction maps $B(G)\to B(\calF)$ and $T(\calO G)\to B(\calF)$ are split surjective.

\bigskip\noindent
{\bf Theorem B}\quad{\it
Let $S$ be a Sylow $p$-subgroup of $G$ and let $\calF:=\calF_S(G)$ denote the resulting fusion system on $S$. There exists a ring homomorphism $t_S^G\colon B(\calF)\to B(G)$ such that the compositions
\begin{equation}\label{eqn splitting 1}
  \res^G_S\circ t_S^G\colon B(\calF)\to B(G) \to B(\calF)
\end{equation}
and
\begin{equation}\label{eqn splitting 2}
   \lambda_S^{-1}\circ\res^G_S\circ\lambda_G\circ t_S^G\colon   
   B(\calF)\to B(G)\to T(\calO G)\to T(\calO S) \myiso B(S)
\end{equation}
are equal to the identity map on $B(\calF)$. In particular, one obtains decompositions
\begin{equation*}
   B(G)=\ker(\res^G_S)\oplus t_S^G(B(\calF))\quad\text{and} \quad 
   T(\calO G)=\ker(\lambda_S^{-1}\circ\res^G_S)\oplus (\lambda_G\circ t_S^G)(B(\calF))
\end{equation*}
into ideals $\ker(\res^G_S)$ of $B(G)$ and $\ker(\lambda_S^{-1}\circ\res^G_S)$ of $T(\calO G)$ and subrings isomorphic to $B(\calF)$.
}

\bigskip
Taking the $\calO$-dual $M^\circ:=\Hom_\calO(M,\calO)$ of a $p$-permutation $\calO G$-module induces a ring automorphism $-^\circ\colon T(\calO G)\to T(\calO G)$, $[M]\mapsto [M^\circ]$, of order $2$. We define the {\em orthogonal unit group} of $T(\calO G)$ as
\begin{equation*}
   \mathrm{O}(T(\calO G)):=\{u\in T(\calO G)\mid u\cdot u^\circ=1\}\,.
\end{equation*}
Thus, $\mathrm{O}(T(\calO G))$ is a subgroup of the unit group of the commutative ring $T(\calO G)$. It turns out that $\mathrm{O}(T(\calO G))$ consists precisely of the units of $T(\calO G)$ of finite order, see Remark~\ref{rem orthogonal units}(b). Orthogonal units of $T(\calO G)$ yield $p$-permutation autoequivalences of $\calO G$, see \cite{BoltjePerepelitsky2020} and see Remark~\ref{rem orthogonal units}(a), and are therefore of particular interest.

\bigskip\noindent
{\bf Theorem C}\quad {\it
Let $S$ be a Sylow $p$-subgroup of $G$ and let $\calF:=\calF_S(G)$ be the associated fusion system on $S$.
One has a direct product decomposition
\begin{equation}\label{eqn orthogonal unit decomposition}
      \mathrm{O}(T(\calO G)) \cong B(\calF)^\times \times \Bigl(\prod_{P\in\scrS_p(G)} \Hom(N_G(P)/P,F^\times)\Bigr)'
\end{equation}
where the second factor is defined as the set of all tuples
\begin{equation*}
   (\varphi_P)\in\Bigl(\prod_{P\in\scrS_p(G)} \Hom(N_G(P)/P,F^\times)\Bigr)^G
\end{equation*}
satisfying
\begin{equation}\label{eqn varphi coherence condition}
   \varphi_P(xP)=\varphi_{P\langle x_p\rangle}(xP\langle x_p\rangle)
\end{equation}
for all $P\in\scrS_p(G)$ and $x\in N_G(P)$.
}

\bigskip
Theorem C is proved in Section~\ref{sec orthogonal units}. It uses Theorems A and B.
In Proposition~\ref{prop SO decomposition} we prove a further direct product decomposition of the second factor in
(\ref{eqn orthogonal unit decomposition}). It is well-known that if $p$ is odd then $B(S)^\times=\{\pm1\}$ and consequently $B(\calF)^\times=\{\pm 1\}$. This follows for instance from \cite[Proposition~6.5]{Yoshida1990}. See also \cite[Lemma~11.2.5]{Bouc2010} for a direct elementary proof. If $p=2$, $B(S)^\times$ has been described explicitly by Bouc (see \cite[Section 11.2]{Bouc2010}) and $B(\calF)^\times$ has been described in some situations by Barsotti and Carman in \cite{BarsottiCarman2020}.

\medskip
In Section~\ref{sec Yoshida criterion} we prove a criterion for a potential tuple of species values to come from an orthogonal unit of $T(\calO G)$. A similar result was proved by Yoshida in \cite[Proposition~6.5]{Yoshida1990} for units of the Burnside ring. In Section~\ref{sec examples} we give more explicit descriptions of $\mathrm{O}(T(\calO G))$ for some special cases of $G$.

\bigskip
{\bf Notation}\quad Apart from the notations already introduced, we will use the following standard notations:

\smallskip
For any ring $A$, we denote by $A^\times$ the unit group of $A$, by $Z(A)$ the center of $A$, and by $\lmod{A}$ the category of finitely generated left $A$-modules. For $M,N\in\lmod{A}$ we write $M\mid N$ if $M$ is isomorphic to a direct summand of $N$.

\smallskip
For any positive integer $n$ we denote by $n_p$ the largest power of $p$ dividing $n$ and set $n_{p'}=n/n_p$. For any group element $x$ we denote by $x_p$ its $p$-part and by $x_{p'}$ its $p'$-part. 

\smallskip
By $a\mapsto\bar{a}$ we denote the natural epimorphisms $\calO\to F$ and $\calO G\to FG$.

\smallskip
We set $\lexp{g}{a}:=gag^{-1}$ for any $g\in G$ and $a\in G$ or $a\in kG$ for a commutative ring $k$. 
For subgroups $H_1$ and $H_2$ of $G$ we write $H_1\le_G H_2$ if $H_1\le \lexp{g}{H_2}$ for some $g\in G$.

\smallskip
For any group $H$, any commutative ring $k$ and any $\varphi\colon H\to k^\times$, we denote by $k_\varphi$,the $kH$-module  whose underlying $k$-module is equal to $k$ and on which $H$ acts via multiplication by $\varphi(h)$. If $\varphi$ is the trivial homomorphism we also write $k_H$ for the trivial $kH$-module.


\section{Preliminaries}

\begin{nothing}\label{noth prel}
(a) Recall that a {\em $p$-permutation} $\calO G$-module is a finitely generated $\calO G$-module $M$ with the property that $\Res^G_P(M)$ is a permutation $\calO P$-module for every $p$-subgroup $P$ of $G$. Equivalently, every indecomposable direct summand of $M$ has trivial source. Another equivalent condition is that $M$ is isomorphic to a direct summand of a finitely generated permutation $\calO G$-module. The class of $p$-permutation modules over $\calO$ is closed under restriction, induction, $\calO$-duals, tensor products $-\otimes_\calO -$, direct sums, and taking direct summands and Green correspondents.

\smallskip
(b) Let $e\in Z(\calO G)$ be an idempotent. We denote the category of $p$-permutation $\calO Ge$-modules by $\ltriv{\calO Ge}$ and denote by $T(\calO Ge)$ the Grothendieck group of the monoid of isomorphism classes of $p$-permutation $\calO Ge$-modules with respect to the direct sum operation. By the Krull-Schmidt Theorem, $T(\calO Ge)$ has a standard $\ZZ$-basis consisting of elements $[M]$, where $M$ runs through a set of representatives of the isomorphism classes of indecomposable $\calO Ge$-modules with trivial source. When $e=1$, $T(\calO G)$ is a commutative ring with multiplication induced by $-\otimes_\calO-$ and it is also known as the {\em trivial source ring} of $G$ over $\calO$. We always consider $T(\calO Ge)$ as a subgroup of $T(\calO G)$.

\smallskip
(c) Similarly, one defines $p$-permutation $FG\ebar$-modules and the associated trivial source group $T(FG\ebar)$. The above statements hold as well over $F$. The functor $F\otimes_\calO -$ induces a ring isomorphism $T(\calO G)\to T(FG)$ which sends the standard basis of $T(\calO Ge)$ bijectively onto the standard basis of $T(FG\ebar)$ and preserves vertices (see \cite{Broue1985}). By $Pr(\calO G e)$  we denote the subgroup of $T(\calO Ge)$ generated by the standard basis elements $[M]$, where $M$ is a projective indecomposable $\calO Ge$-module. Similarly we define $Pr(FG\ebar)$.

\smallskip
(d) Note that one has a commutative diagram
\begin{diagram}[70]
  \movevertex(-20,0){Pr(\calO Ge)} & \movevertex(-10,0){\quad \subseteq\quad}  &   T(\calO Ge) & 
               \movevertex(10,0){\Ear[40]{\kappa_G}} & \movevertex(20,0){R(K Ge)} &&
   \movevertex(-10,0){\saR{\wr}} & & \saR{\wr} & & \movevertex(20,0){\saR{d_G}} &&
   \movevertex(-20,0){Pr(FG\ebar)} & \movevertex(-10,0){\quad \subseteq \quad} &T(FG\ebar) & 
               \movevertex(10,0){\Ear[40]{\eta_G}} & \movevertex(20,0){R(FG\ebar)\,,} &&
\end{diagram}
where $R(KGe)$ denotes the group of virtual characters of $KGe$-modules, $R(FGe)$ denotes the group of virtual Brauer characters of $FGe$-modules, $\kappa_G$ is induced by the functor $K\otimes_\calO -\colon \ltriv{\calO Ge}\to\lmod{KGe}$, $\eta_G$ is the canonical map, the left two vertical maps are induced by the functor $F\otimes_\calO-$, and the right vertical map is the decomposition map. By \cite[Th\'eor\`emes~34 et 36]{Serre1978}, the map $\kappa_G\colon Pr(\calO Ge)\to R(KGe)$ is injective and its image consists precisely of those virtual characters in $R(KGe)$ which vanish on $p$-singular elements.
\end{nothing}

\begin{lemma}\label{lem Brauer construction 1}
Suppose that $P$ is a normal $p$-subgroup of a finite group $H$. Let $M$ be an indecomposable $p$-permutation $\calO H$-module with vertex $Q$. The following are equivalent:

\smallskip
{\rm (i)} $P$ acts trivially on $M$, i.e., $\Res^H_P(M)\cong \calO_P\oplus \cdots\oplus \calO_P$.

\smallskip
{\rm (ii)} $\calO_P\mid \Res^H_P(M)$.

\smallskip
{\rm (iii)} $P\le Q$.
\end{lemma}

\begin{proof}
Obviously, (i) implies (ii). If (ii) holds, then $\calO_P\mid \Res^H_P(\Ind_Q^H(\calO_Q))$. The Mackey formula and the fact that $\calO_P$ has vertex $P$ now implies that (iii) holds. Finally, if (iii) holds then $\Res^H_P(M)\mid
\Res^H_P(\Ind_Q^H(\calO_Q))$ and the Mackey formula imply (i).
\end{proof}

An analogous version of Lemma~\ref{lem Brauer construction 1} with identical proof holds over $F$ instead of $\calO$.

\begin{nothing}\label{noth prime construction}
Let $P\in\scrS_p(G)$. Recall from \cite{Broue1985} that the {\em Brauer construction} with respect to $P$ is a functor
\begin{equation*}
   -(P)\colon\ltriv{FG}\to\ltriv{F[N_G(P)/P]}\,,\quad M\mapsto M(P)\,.
\end{equation*}
If $M\in\ltriv{FG}$ is indecomoposable and $P$ is a vertex of $M$ then the inflated module $\Inf_{N_G(P)/P}^{N_G(P)}( M(P))\in\ltriv{N_G(P)}$ is the Green correspondent of $M$. Moreover, for $Q\in \scrS_p(G)$, one has $M(Q)\neq \{0\}$ if and only if $Q\le_G P$.

The commutative diagram
\begin{diagram}[75]
   \movevertex(-40,0){T(\calO G)} & \movearrow(-32,0){\Ear[40]{-(P)}} & T(\calO[N_G(P)/P]) &&
   \movearrow(-40,0){\Sar{\wr}} & & \saR{\wr} &&
   \movevertex(-40,0){T(FG)} & \movearrow(-32,0){\Ear[40]{-(P)}} & T(F[N_G(P)/P]) &&
\end{diagram}
defines  a unique ring homomorphism $T(\calO G)\to T(\calO[N_G(P)/P])$ which we abusively denote again by $-(P)$. Note however, that this map does not come from a functor $\ltriv{\calO G}\to\ltriv{\calO[N_G(P)/P]}$. For $M\in\ltriv{\calO G}$ the image $[M](P)$ of $[M]\in T(\calO G)$, can be constructed as follows: Decompose
\begin{equation*}
   \Res^G_{N_G(P)}(M)=V_1\oplus\ldots\oplus V_n
\end{equation*}
into indecomposable $\calO{N_G(P)}$-modules $V_1,\ldots,V_n$ and let $I$ denote the set of all $i\in\{1,\ldots,n\}$ such that $V_i$ satisfies the equivalent conditions of Lemma~\ref{lem Brauer construction 1} (with $H=N_G(P)$), i.e, $P$ acts trivially on $V_i$. For each $i\in I$ one has $V_i\cong\Inf_{N_G(P)/P}^{N_G(P)}(\Vtilde_i)$ for an indecomposable $\Vtilde_i\in\ltriv{\calO[N_G(P)/P]}$, which is uniquely determined up to isomorphism by $V_i$. Then $[M](P)=\sum_{i\in I} [\Vtilde_i]$. By abuse of notation we also set
\begin{equation*}
   M(P):=\bigoplus_{i\in I} \Vtilde_i\in\ltriv{\calO[N_G(P)/P]}\,,
\end{equation*}
which, by the Krull-Schmidt Theorem is uniquely determined up to isomorphism.
\end{nothing}

\begin{lemma}\label{lem Brauer construction 2}
Let $M$ be a $p$-permutation $\calO G$-module and let $P\in\scrS_p(G)$.
For any subgroup $H$ satisfying $P\le H\le N_G(P)$ one has 
\begin{equation*}
   \Res^{N_G(P)/P}_{H/P} (M(P)) \cong \bigl(\Res^G_H(M)\bigr)(P)
\end{equation*}
as $\calO[H/P]$-modules.
\end{lemma}

\begin{proof}
Let $V$ be an indecomposable direct summand of $\Res^G_{N_G(P)}(M)$ and let $W$ be an indecomposable direct summand of $\Res^{N_G(P)}_H(V)$. By Lemma~\ref{lem Brauer construction 1}, the $\calO[N_G(P)]$-module $V$ satisfies the conditions in Lemma~\ref{lem Brauer construction 1} if and only if the $\calO H$-module $W$ does. In fact, if $V$ satisfies $(i)$ then $W$ does, and if $W$ satisfies (ii) then $V$ does. The result follows.
\end{proof}


\section{Proof of Theorem~A}\label{sec Thm A}

In this section we will prove Theorem~A and generalize it to statements that involve sums of blocks and restricted vertices.

\smallskip
Using the construction from \ref{noth prime construction}, we first define the homomorphism $\beta_G$.
Note that $G$ acts on the ring $\prod_{P\in\scrS_p(G)} R(K[N_G(P)/P])$ via ring automorphisms by
\begin{equation*}
   \lexp{g}{((\chi_P)_{P\in\scrS_p(G)})}:=\bigl(\lexp{g}{(\chi_{\lexp{g^{-1}}{P}}})\bigr)_{P\in\scrS_p(G)}\,.
\end{equation*}
Here, $\lexp{g}{(\chi_{\lexp{g^{-1}}{P}}})\in R(K[N_G(P)/P])$ denotes the character arising from $\chi_{\lexp{g^{-1}}{P}}$ via restriction along the isomorphism $N_G(P)/P\to N_G(\lexp{g^{-1}}{P})/\lexp{g^{-1}}{P}$ induced by $x\mapsto g^{-1}xg$ for $x\in N_G(P)$.
Thus, a tuple $(\chi_P)$ is a $G$-fixed point if and only if $\lexp{g}{\chi_P}=\chi_{\lexp{g}{P}}$ for all $g\in G$ and $P\in\scrS_p(G)$. The homomorphism $\beta_G$ is now defined as
\begin{equation*}
   \beta_G\colon T(\calO G)\to \Bigl(\prod_{P\in\scrS_p(G)} R(K[N_G(P)/P])\Bigr)^G\,,\quad 
   \omega\mapsto \bigl(\kappa_{N_G(P)/P}({\omega(P)})\bigr)_{P\in\scrS_p(G)}\,.
\end{equation*}
Since both the maps $-(P)\colon T(\calO G) \to T(\calO[N_G(P)/P])$ and $\kappa_{N_G(P)/P}$ commute with $G$-conjugations, the image of $\beta_G$ is contained in the $G$-fixed points, and since both maps are ring homomorphisms, $\beta_G$ is a ring homomorphism.

\bigskip
\begin{proof}{\em of Theorem A.}\quad First note that $\beta_G$ is injective, since, by the commutative diagram in \ref{noth prel}(d), the diagram
\begin{diagram}
   \movevertex(-120,0){T(\calO G)} & \movearrow(-110,0){\Ear[50]{(-(P))}} &  
        \movevertex(-50,0){\Bigl(\mathop{\prod}\limits_{P\in\scrS_p(G)} T(\calO[N_G(P)/P])\Bigr)^G} &
        \movearrow(10,0){\Ear[50]{(\kappa_{N_G(P)/P})}} &
        \movevertex(75,0){\Bigl(\mathop{\prod}\limits_{P\in\scrS_p(G)} R(K[N_G(P)/P])\Bigr)^G} &&
   \movearrow(-120,0){\Sar{\wr}} & & \movearrow (-50,0){\Sar{\wr}} & & \movearrow(75,0){\saR{(d_{N_G(P)/P})}} &&
   \movevertex(-120,0){T(F G)} & \movearrow(-110,0){\Ear[50]{(-(P))}} &  
        \movevertex(-50,0){\Bigl(\mathop{\prod}\limits_{P\in\scrS_p(G)} T(F[N_G(P)/P])\Bigr)^G} &
        \movearrow(10,0){\Ear[50]{(\eta_{N_G(P)/P})} }&
            \movevertex(75,0){\Bigl(\mathop{\prod}\limits_{P\in\scrS_p(G)} R(F[N_G(P)/P])\Bigr)^G} &&
\end{diagram}
is commutative and the composition of the two bottom maps is injective by a theorem of Conlon, see \cite[Proposition~5.5.4]{Benson1998}. 

\smallskip
Next we show that every element in the image of $\beta_G$ satisfies the condition (C) in Theorem~A. This can be derived from \cite[Lemma~6.2]{Rickard1996}, but we give a detailed proof for the reader's convenience. Let $M\in\ltriv{\calO G}$ and $(\chi_P)=\beta_G([M])$. Let $P\in\scrS_p(G)$, and $x\in N_G(P)$. We need to show that $\chi_P(xP) = \chi_Q(xQ)$  with $Q:=P\langle x_p \rangle$. By Lemma~\ref{lem Brauer construction 2} with $H:=P\langle x \rangle$, we may assume that $G=P\langle x \rangle$. We may also assume that $M$ is indecomposable and $P<Q$. Let $R$ be a vertex of $M$. If $R=Q$ then $M(P)=M=M(Q)$ and therefore $\chi_P(xP)=\chi_Q(xQ)$. If $R\neq Q$ then $R<Q$ and therefore $M(Q)=\{0\}$ and $\chi_Q(xQ)=0$. In order to show that $\chi_P(xP)=0$, we distinguish the two cases $P\le R$ and $P\not\le R$. In the latter case, $M(P)=\{0\}$ and we are done. So assume that $P\le R<Q$. Then, with $\calOtilde$ denoting a finite extension of $\calO$ containing a root of unity of order $|\langle x_{p'}\rangle|$, we have
\begin{equation*}
   \calOtilde\otimes_{\calO}M\mid \Ind_R^G(\calOtilde) \cong \Ind_{R\langle x_{p'} \rangle}^G \Ind_R^{R \langle x_{p'} \rangle} ( \calOtilde ) \cong
   \bigoplus_{\varphi} \Ind_{R \langle x_{p'} \rangle}^G (\calOtilde_\varphi)\,,
\end{equation*}
where $\varphi$ runs through all homomorphisms $R\langle x_{p'} \rangle\to\calOtilde^\times$ with $\varphi|_R=1$. Green's indecomposablility theorem implies that $\calOtilde\otimes_\calO M$ is isomorphic to a direct sum of modules of the form $\Ind_{R \langle x_{p'} \rangle}^G (\calOtilde_\varphi)$, for homomorphisms $\varphi$ as above. Therefore, $M(P)=M$ and $\chi_P(xP)$ equals a sum of characters of $\calOtilde G$-modules of the form $\Ind_{R\langle x_{p'} \rangle}^G(\calOtilde_\varphi)$ evaluated at $x$. The latter evaluation equals $0$, since $x$ is not $G$-conjugate to an element of $R \langle x_{p'} \rangle$ noting that $R<Q$.

\smallskip
Conversely, assume that we have a family of virtual characters $\chi_P\in R(K[N_G(P)/P])$ satisfying 
\begin{equation}\label{eqn prime condition}
   \lexp{g}{\chi_P}= \chi_{\lexp{g}{P}}\quad\text{and}\quad 
   \chi_P(xP) = \chi_{P\langle x_p\rangle}(xP\langle x_p\rangle)
\end{equation}
for all $P\in\scrS_p(G)$, $g\in G$, and $x\in N_G(P)$. We will construct an element $\omega\in T(\calO G)$ with $\beta_G(\omega) = (\chi_P)_{P\in\scrS_p(G)}$ inductively. Let $P_1,\ldots,P_n\in\scrS_p(G)$ be representatives of the $G$-conjugacy classes and assume that $|P_1|\ge|P_2|\ge\cdots\ge|P_n|$. We will show by induction on $i\in\{1,\ldots, n\}$ that there exists an element $\omega_i\in T(\calO G)$ with the following two properties:

\smallskip
(i) $\omega_i$ is a $\ZZ$-linear combination of the standard basis elements $[M]$ with vertex $P_i$. 

\smallskip
(ii) The $P_i$-component of $\beta_G(\omega_1 +\cdots+ \omega_i)$ is equal to $\chi_{P_i}$.

\smallskip
$i=1$: Note that $P:=P_1$ is a Sylow $p$-subgroup of $G$. Since $N_G(P)/P$ is a $p'$-group, we have $Pr(\calO[N_G(P)/P])=T(\calO[N_G(P)/P])$ and $\kappa_{N_G(P)/P}\colon Pr(\calO[N_G(P)/P])\to R_K(N_G(P)/P)$ is an isomorphism. Thus, there exists a $\ZZ$-linear combination $\gamma:=a_1[V_1]+\cdots+a_k[V_k]$ with indecomposable projective $\calO[N_G(P)/P]$-modules $V_1,\ldots, V_k$ such that $\kappa_{N_G(P)/P}(\gamma)= \chi_{P}$. Then each $W_r:=\Inf_{N_G(P)/P}^{N_G(P)}(V_r)$, $r=1,\ldots,k$, is an indecomposable $p$-permutation $\calO[N_G(P)]$-module with vertex $P$. For $r=1,\ldots,k$, let $M_r\in\ltriv{OG}$ be the Green correspondent of $W_r$ and set $\omega_1:=a_1[M_1]+\cdots+a_k[M_k]\in T(\calO G)$. Then $\omega_1$ satisfies (i) and (ii).

\smallskip
$i\to i+1$: Set $P:=P_{i+1}$, let $x\in N_G(P)$ with $x_p\notin P$, and set $Q:=P\langle x_p\rangle$. Then $Q>P$. Let $\beta_G(\omega_1+\cdots+\omega_i)_P\in R_K(N_G(P)/P)$ denote the $P$-component of $\beta_G(\omega_1+\cdots+\omega_i)$.

{\em Claim:} $\beta_G(\omega_1+\cdots+\omega_i)_P(xP)=\chi_P(xP)$. To prove the claim, let $j\in\{1,\ldots,i\}$ with $Q=_G P_j$. We first show that 
\begin{equation}\label{eqn j+1}
   \beta_G(\omega_{j+1}+\cdots+\omega_i)_P (xP) = 0\,.
\end{equation}
If $\omega_{j+1}+\cdots+\omega_i= 0$ then this is obvious. So assume that $M$ is an indecomposable $p$-permutation $\calO G$-module such that $[M]$ appears in $\omega_{j+1}+\cdots+\omega_i$ with nonzero coefficient and let $W$ be an indecomposable $\calO N_G(P)$-module that appears in $M(P)$. It suffices to show that the $\chi_W(x)=0$, where $\chi_W$ denotes the character of $W$. So assume that $\chi_{W}(x)\neq 0$, which implies $x_p\in R$ for a vertex $R$ of $W$. Since $W\mid M(P)\mid \Res^G_{N_G(P)}(M)$, we have $R\ge P$. Thus $Q=P\langle x_p\rangle \le R$. But by (i), $M$ has a vertex $R'$ with $|R'|\le|P_j|=|Q|$ and since $W\mid M(P)\mid \Res^G_{N_G(P)}(M)$, we also have $R\le_G R'$. Thus, $|R|\le|Q|$ so that $Q\le R$ implies $R=Q$. Thus, $M$ has vertex $P_j=_GQ$, in contradiction to $M$ appearing in $\omega_{j+1}+\cdots+\omega_i$. Hence, (\ref{eqn j+1}) holds. Next note that since $\beta_G(\omega_1+\cdots+\omega_j)$ satisfies the coherence condition (C), we obtain
\begin{equation}\label{eqn j}
   \beta_G(\omega_1+\cdots+\omega_j)_P(xP) = \beta_G(\omega_1+\cdots+\omega_j)_Q(xQ)\,.
\end{equation}
Further, by the induction hypothesis for $j\le i$, and since the tuple $(\chi_P)$ satisfies the coherence condition we have
\begin{equation}\label{eqn P Q}
   \beta_G(\omega_1+\cdots+\omega_j)_Q(xQ)=\chi_Q(xQ) = \chi_P(xP)\,.
\end{equation}
Combining the equations (\ref{eqn j+1}), (\ref{eqn j}) and (\ref{eqn P Q}) now proves the claim.

Next we set
\begin{equation*}
   \psi:=\chi_P-\beta_G(\omega_1+\cdots+\omega_i)_P\in R_K(N_G(P)/P)\,,
\end{equation*}
which is a virtual character that vanishes on $p$-singular elements by the claim. By \cite[Th\'eor\`emes 34 et 36]{Serre1978}, see also \ref{noth prel}(d), we can write $\psi=\kappa_{N_G(P)/P}(\gamma)$ for a unique $\ZZ$-linear combination $\gamma=a_1[V_1]+\cdots + a_k[V_k]$ with indecomposable projective $\calO[N_G(P)/P]$-modules $V_1,\ldots, V_k$. For $r=1,\ldots,k$, set $W_r:=\Inf_{N_G(P)/P}^{N_G(P)}(V_r)$. Then each $W_r$ has vertex $P$. For each $r=1,\ldots,k$, let $M_r$ denote the Green correspondent of $W_r$. Then each $M_r$ is an indecomposable $p$-permutation $\calO G$-module with vertex $P$. We set 
\begin{equation*}
   \omega_{i+1}:=a_1[M_1]+\cdots+a_k[M_k]\,.
\end{equation*}
Thus, $\omega_{i+1}$ satisfies (i). Moreover, it satisfies (ii), since 
\begin{equation*}
   \beta_G(\omega_{i+1})_P = \kappa_{N_G(P)/P} (a_1[V_1]+\cdots a_k[V_k]) = \psi
\end{equation*}
and therefore
\begin{equation*}
   \beta_G(\omega_1+\cdots+\omega_{i+1})_P
   = \beta_G(\omega_1+\cdots+\omega_i)_P+ \psi = \chi_P\,.
\end{equation*}
This completes the induction proof.

\smallskip
Finally, set $\omega:=\omega_1+\cdots+\omega_n$ and let $i\in\{1,\ldots,n\}$. Then we have
\begin{equation*}
   \beta_G(\omega)_{P_i} = \beta_G(\omega_1+\cdots+\omega_i)_{P_i}=\chi_{P_i}\,.
\end{equation*}
In fact, the first equation follows from $\omega_j(P_i)=0$ for all $j\in\{i+1,\ldots,n\}$, since $P_i\not\le_G P_j$, see \ref{noth prime construction}, and the second equation follows from (ii). Finally, since the family $(\chi_P)$ is $G$-fixed we also have $\beta_G(\omega)_P=\chi_P$ for any $P\in\scrS_p(G)$.
\end{proof}

\begin{notation}\label{not prime}
We will denote from now on by
\begin{equation}
   \Bigl(\prod_{P\in\scrS_p(G)} R\bigl(K[N_G(P)/P]\bigr)\Bigr)'
\end{equation}
the set of all tuples $(\chi_P)$ with $\chi_P\in R(K[N_G(P)/P])$, for $P\in\scrS_p(G)$, which satisfy the two conditions in (\ref{eqn prime condition}).
\end{notation}

\begin{nothing}
Next we derive a block-wise version of Theorem A. Let $e\in Z(\calO G)$ be an idempotent. For $P\in\scrS_p(G)$ denote by $e(P)$ the unique idempotent in $Z(\calO[N_G(P)])$ with $\br_P(e)=\overline{e(P)}$, where $\br_P\colon Z(\calO G)\to Z(F[N_G(P)])$ denotes the Brauer homomorphism. It is a well-known property of the Brauer construction that 
\begin{equation}\label{eqn Brauer construction and blocks}
   e(P)M(P)=M(P)\ \text{for every $M\in\ltriv{\calO Ge}$}\,.
\end{equation}
Note that, since $\lexp{g}{e(P)}=e(\lexp{g}{P})$ for all $g\in G$ and $P\in\scrS_p(G)$, the subgroup $\prod_{P\in\scrS_p(G)} R\bigl(K[N_G(P)/P]e(P)\bigr)$ of $\prod_{P\in\scrS_p(G)}R\bigl(K[N_G(P)/P]\bigr)$ is $G$-stable.
Therefore, if we denote by
\begin{equation*}
   \Bigl(\prod_{P\in\scrS_p(G)} R\bigl(K[N_G(P)/P]e(P)\bigr)\Bigr)'
\end{equation*}
the set of all tuples $(\chi_P)$ with $\chi_P\in R(K[N_G(P)/P]e(P))$, for $P\in\scrS_p(G)$, which satisfy the two conditions in (\ref{eqn prime condition}) then $\beta_G$ restricts to an injective group homomorphism
\begin{equation}\label{eqn blockwise beta}
   \beta_G\colon T(\calO Ge)\to \Bigl(\prod_{P\in\scrS_p(G)} R\bigl(K[N_G(P)/P]e(P)\bigr)\Bigr)'\,,
\end{equation}
where we consider $e(P)$ also as a central idempotent of $K[N_G(P)/P]$ via the canonical $K$-algebra epimorphism $K[N_G(P)]\to K[N_G(P)/P]$. Note that $T(\calO Ge)$ is in general not a subring of $T(\calO G)$.
\end{nothing}

\begin{corollary}\label{cor Thm A blockwise}
For every idempotent $e\in Z(\calO G)$, the map in (\ref{eqn blockwise beta}) is an isomorphism.
\end{corollary}

\begin{proof}
It suffices to prove the surjectivity. So let $(\chi_P)$ be a tuple in the codomain of (\ref{eqn blockwise beta}). Then, by Theorem~A there exists $\omega\in T(\calO G)$ with $\beta_G(\omega)=(\chi_P)$. The decomposition $T(\calO G) = T(\calO Ge)\oplus T(\calO G(1-e))$ allows us to write $\omega=\omega e + \omega(1-e)$ with $\omega e\in T(\calO Ge)$ and $\omega(1-e)\in T(\calO G(1-e))$. The additivity of $\beta_G$ and the property (\ref{eqn Brauer construction and blocks}) for the idempotents $e$ and $1-e$ now imply that $\beta_G(\omega(1-e))=0$ and $\beta(\omega e)=(\chi_P)$.
\end{proof}

\begin{notation}
Let $\scrX$ be a subset of $\scrS_p(G)$ that is closed under $G$-conjugation and under taking subgroups and let $e\in Z(\calO G)$ be an idempotent. We denote by $T(\calO Ge\,|\, \scrX)$ the $\ZZ$-span of the standard basis elements $[M]$, where $M\in\ltriv{\calO Ge}$ is indecomposable and has a vertex belonging to $\scrX$. Moreover, we denote by
\begin{equation*}
   \Bigl(\prod_{P\in\scrX} R\bigl(K[N_G(P)/P]e(P)\bigr)\Bigr)'
\end{equation*}
the set of all tuples $(\chi_P)$ with $\chi_P\in R(K[N_G(P)/P]e(P))$, for $P\in \scrX$, satisfying
\begin{equation*}
      \lexp{g}{\chi_P}= \chi_{\lexp{g}{P}}\quad\text{and}\quad 
   \chi_P(xP) = \chi_Q(xQ)
\end{equation*}
for all $g\in G$, $P\in\scrX$ and $x\in N_G(P)$, with $Q:=P\langle x_p \rangle$, where we interpret $\chi_Q$ as $0$ if $Q\notin\scrX$. Note that for $M\in T(\calO Ge\,|\, \scrX)$, one has $M(P)=\{0\}$ whenever $P\in\scrS_p(G)\smallsetminus\scrX$. Thus, the map $\beta_G$ restricts to an injective group homomorphism
\begin{equation}\label{eqn blockwise beta with X}
   \beta_G\colon T(\calO Ge\,|\, \scrX) \to \Bigl(\prod_{P\in\scrX} R\bigl(K[N_G(P)/P]e(P)\bigr)\Bigr)'\,,
\end{equation}
whose codomain we consider as a subset of $\prod_{P\in\scrS_p(G)} R(K[N_G(P)/P]e(P))$ by adding $0$'s in the missing components.
\end{notation}

\begin{corollary}\label{cor Thm A with scrX}
Let $\scrX\subseteq\scrS_p(G)$ be a subset that is closed under taking $G$-conjugates and under taking subgroups and let $e\in Z(\calO G)$ be an idempotent. Then the map in (\ref{eqn blockwise beta with X}) is an isomorphism.
\end{corollary}

\begin{proof}
It suffices to show the surjectivity. Suppose that $(\chi_P)$ belongs to the codomain of the map (\ref{eqn blockwise beta with X}), and view this tuple as element in $\prod_{P\in\scrS_p(G)} R(K[N_G(P)/P]e(P))$ by adding $0$'s in the components indexed by $P\in\scrS_p(G)\smallsetminus\scrX$. Then $(\chi_P)$ belongs to the codomain of the map in (\ref{eqn blockwise beta}). 
By Corollary~\ref{cor Thm A blockwise} there exists $\omega\in T(\calO Ge)$ with $\beta_G(\omega)=(\chi_P)$ for all $P\in\scrS_p(G)$. We'll show that $\omega\in T(\calO G\, | \,\scrX)$. 
Let $\scrY$ denote the set of all $P\in\scrS_p(G)$ such that there exists an indecomposable $p$-permutation $\calO G$-module $M$ having a vertex that contains $P$ and occurring with nonzero coefficient in $\omega$. Then $\scrY$ is closed under $G$-conjugation and taking subgroups. It suffices to show that $\scrY\subseteq\scrX$. Let $P$ be a maximal element in $\scrY$. 
It suffices to show that $P\in\scrX$. Suppose this is not the case. Then $0=\chi_P=\kappa_{N_G(P)/P}(\omega(P))$. Moreover, $\omega(P)\in Pr(\calO[N_G(P)/P])$, since $P$ was chosen maximal in $\scrY$. Since $\kappa_{N_G(P)/P}$ is injective on $Pr(\calO[N_G(P)/P])$, this implies that $\omega(P)=0$. 
We can write $\omega=\omega_1+\omega_2$, such that $P$ is a vertex of the indecomposable modules appearing in $\omega_1$ and $P$ is not a vertex of the indecomposable modules appearing in $\omega_2$. By the maximality of $P$ in $\scrY$, we obtain $\omega_1\neq 0$ but $\omega_1(P)=\omega(P)=0$.
We write $\omega_1=a_1[M_1]+\cdots+a_k[M_k]$ with representatives $M_1,\ldots M_k$ of the isomorphism classes of indecomposable $p$-permutation $\calO Ge$-modules with vertex $P$ and $a_1,\ldots,a_k\in\ZZ$, we obtain $\omega_1(P)=a_1[V_1]+\cdots+a_k[V_k]$, where each $V_r$ is an indecomposable projective $\calO[N_G(P)/P]e(P)$-module whose inflation $W_r$ to $N_G(P)$ is the Green correspondent of $M_r$, for $r=1,\ldots,k$. Since the Green correspondence is bijective and the inflation map $T(\calO[N_G(P)/P]e(P))\to T(\calO[N_G(P)]e(P))$ is injective, $V_1,\ldots,V_r$ are pairwise non-isomorphic. Thus, $a_1=\cdots=a_k=0$ and therefore $\omega_1=0$, a contradiction.
\end{proof}

\begin{remark}
(a) It is worth mentioning that in Theorem A, the field $K$ can be any finite extension of $\QQ_p$, including $\QQ_p$ itself.

\smallskip
(b) Suppose that $K$ contains a root of unity of order $\exp(G)$, the exponent of $G$, so that $R(KG)$ is isomorphic to the character ring of $G$. The image of $\beta_G$ followed by the projection map $\pi_1$ onto $R(KG)$ is equal to $\kappa_G(T(\calO G))$. The coherence condition (C) in Theorem A implies the well-known fact that characters of $p$-permutation $\calO G$-modules are $p$-rational, i.e., with values contained in the group $\QQ_p(\zeta)$, where $\zeta$ is a root of unity of $G$ of order $\exp(G)_{p'}$. The latter is also the group $R(\QQ_p(\zeta)G)$ of characters of $\QQ_p(\zeta)$-modules by a Theorem of Brauer, see \cite[Lemma~$1^*$]{Brauer1941}. It is shown in \cite{McHugh2021} that the factor group $R(\QQ_p(\zeta)G)/\kappa_G(T(\calO G))$ is an elementary abelian $2$-group (trivial if $p$ is odd) and that these groups (for varying $G$) form a uniserial fibered biset functor over the field $\FF_2$ whose composition factors are indexed by the quaternion groups.
\end{remark}


\section{The orthogonal unit group of $T(\calO G)$}\label{sec orthogonal units}

Recall from Section~1 that the orthogonal unit group of $T(\calO G)$ is defined by $\mathrm{O}(T(\calO G))=\{u\in T(\calO G)\mid u\cdot u^\circ=1\}$.

\begin{remark}\label{rem orthogonal units}
(a) Our interest in the orthogonal unit group $\mathrm{O}(T(\calO G))$ stems from the fact that the functor $\Ind_{\Delta(G)}^{G\times G} \colon \ltriv{\calO G}\to \ltriv{\calO[G\times G]}$, where we identify $G$ and the diagonal subgroup $\Delta(G):=\{(g,g)\mid g\in G\}$ of $G\times G$ via $g\mapsto (g,g)$, yields an injective ring homomorphism $\Delta\colon T(\calO G)\to T^\Delta(\calO G,\calO G)$, when we identify left $\calO[G\times G]$-modules with $(\calO G,\calO G)$-bimodules. Here $T^\Delta(\calO G,\calO G)$ denotes the free abelian group on isomorphism classes of indecomposable $(\calO G,\calO G)$-bimodules which are $p$-permutation modules with twisted diagonal vertices when considered as $\calO[G\times G]$-modules. The product on the ring $T^\Delta(\calO G,\calO G)$ is induced by taking the tensor product over $\calO G$. If $u$ is an orthogonal unit of $T(\calO G)$ then its image $\Delta(u)\in T^\Delta(\calO G,\calO G)$ is a $p$-permutation autoequivalence of $\calO G$, see \cite{BoltjePerepelitsky2020}. Thus, $\mathrm{O}(T(\calO G))$ can be considered as a subgroup of the group of $p$-permutation autoequivalences of $\calO G$, a group we are interested in studying in the future.

\smallskip
(b) The group $\mathrm{O}(T(\calO G))$ is actually equal to the torsion subgroup of the unit group $T(\calO G)^\times$, see also \cite[Section~9]{Carman2018}. To see this, consider the map
\begin{equation*}
   \sigma_G\colon T(\calO G)\to \prod_{(P,sP)\in\scrT_p(G)} \ZZ[\zeta]\,,
\end{equation*}
where $\scrT_p(G)$ is the set of all pairs $(P,sP)$ with $P\in\scrS_p(G)$ and $sP$ a $p'$-element of $N_G(P)/P$, and where the $(P,sP)$-component of $\sigma_G(\omega)$ is given by the evaluation of the virtual character $\kappa_{N_G(P)/P}(\omega(P))$ at the element $sP$. Here $\zeta$ can be taken as a root of unity of order $\exp(G)_{p'}$ in some extension of $K$. The map $\sigma_G$ is an injective ring homomorphism (see \cite[Corollary~5.5.5]{Benson1998}) whose components are also called the {\em species} of $T(\calO G)$. It commutes with duals: if the $(P,sP)$-component of $\sigma_G(\omega)$ is equal to $a\in \ZZ[\zeta]$ then the $(P,sP)$-component of $\sigma(\omega^\circ)$ is equal to the \lq complex conjugate\rq\ $\abar$ of $a$, i.e., the image of $a$ under the unique Galois automorphism of $\QQ(\zeta)$ sending $\zeta$ to $\zeta^{-1}$. By \cite[Theorem~4.12]{Washington1982}, every element $a\in\ZZ[\zeta]$ with $a\abar=1$ has finite order. This implies that the image of every element $u\in \mathrm{O}(T(\calO G))$ has finite order, and that also $u$ has finite order by the injectivity of $\sigma_G$. Conversely, if $u\in T(\calO G)^\times$ has finite order then also $\sigma_G(u)$ has finite order. Thus, every component $a_{(P,sP)}$ of $\sigma_G(u)$ is a root of unity and satisfies $a_{(P,sP)}\abar_{(P,sP)}=1$. Since $\sigma_G$ transforms duals into complex conjugates and $\sigma_G$ is injective, $u$ satisfies $u u^\circ=1$.

\smallskip
(c) Taking $K$-duals of $KG$-modules induces a ring automorphism $-^\circ\colon R(KG)\to R(KG)$. Again, we call a unit $u\in R(KG)^\circ$ an {\em orthogonal unit} if $uu^\circ=1$. The orthogonal units of $R(KG)$ form a subgroup $\mathrm{O}(R(KG))$ of the unit group $R(KG)^\times$. Using the species of $R(KG)$, i.e., the injective ring homomorphism $R(KG)\to\prod_{g\in G} \ZZ[\zeta]$ with $\zeta$ a primitive root of unity of order $\exp{G}$ in some extension field of $K$, the same arguments as in (b) show that $\mathrm{O}(R(KG))$ equals the torsion subgroup of $R(KG)^\times$. Moreover, it follows from \cite{Yamauchi1991} that every torsion unit of $R(KG)$ is of the form $\epsilon\cdot \varphi$, where $\epsilon\in\{\pm1\}$ and $\varphi\in\Hom(G,K^\times)$ is a linear character with values in $K$. Therefore, $\mathrm{O}(R(KG))=\{\pm1\}\times \Hom(G,K^\times)$, see also \cite[Section~9]{Carman2018}.
\end{remark}

\begin{proposition}\label{prop beta(u) decomposition}
For each $u\in \mathrm{O}(T(\calO G))$ there exist unique signs $\epsilon_P\in\{\pm1\}$ and homomorphisms $\varphi_P\in\Hom(N_G(P)/P, K^\times)$, for $P\in\scrS_P(G)$, such that
\begin{equation*}
   \beta_G(u)=(\epsilon_P\cdot\varphi_P)_{P\in\scrS_p(G)}\,.
\end{equation*}
\end{proposition}

\begin{proof}
First note that $\beta_G(\omega^\circ)=\beta_G(\omega)^\circ$, for all $\omega\in T(\calO G)$, where we use the notation $(\chi_P)^\circ:=(\chi_P^\circ)$ for any tuple $(\chi_P)\in\prod_{P\in\scrS_p(G)} R(K[N_G(P)/P])$. Therefore, if $u\in \mathrm{O}(T(\calO G))$ then the $P$-component of $\beta_G(u)$ is an orthogonal unit in $R(K[N_G(P)/P])$. The result follows now from Remark~\ref{rem orthogonal units}(c).
\end{proof}

We will continue the study of $\mathrm{O}(T(\calO G))$ after the proof of Theorem~B.

\begin{definition}(\cite{Reeh2015})\label{def B(calF)}
Let $S$ be a Sylow $p$-subgroup of $p$ and let $\calF:=\calF_S(G)$ denote the associated fusion system on $S$. The Burnside ring of $\calF$, denoted by $B(\calF)$ is defined as the subring of the Burnside ring $B(S)$ consisting of all elements $a\in B(S)$ with the property that $\phi_P(a)=\phi_Q(a)$, whenever $P$ and $Q$ are $G$-conjugate subgroups of $S$. Here, for $P\le S$, $\phi_P\colon B(S)\to\ZZ$ is the ring homomorphism defined by $\phi_P([X]):=|X^P|$, whenever $X$ is a finite $S$-set, where $X^P$ denotes the set of $P$-fixed points of $X$.
\end{definition}

\begin{theorem}(\cite[Theorem~7.1]{BarsottiCarman2020})\label{thm BaCa} Let $S$ be a Sylow $p$-subgroup of $G$ and let $\calF:=\calF_S(G)$ denote the associated fusion system on $S$. Then the ring homomorphism $\res^G_S\colon B(G)\to B(\calF)$ induced by restricting a finite $G$-set to $S$, is split surjective. More precisely, there exists a ring homomorphism $t_S^G\colon B(\calF)\to B(G)$ such that $\res^G_S\circ t_S^G=\id_{B(\calF)}$. For $a\in B(\calF)$, the element $t_S^G(a)$ is determined by the equations $\phi_H(t_S^G(a))=\phi_{P_H}(a)$, for all subgroups $H\le G$, where $P_H$ is a subgroup of $S$ which is $G$-conjugate to a Sylow $p$-subgroup of $H$.
\end{theorem}

\begin{proof} {\em of Theorem B.}\quad That (\ref{eqn splitting 1}) is the identity on $B(\calF)$ is already stated in Theorem~\ref{thm BaCa}. And that (\ref{eqn splitting 2}) is the identity on $B(\calF)$ follows from the obvious commutativity of the diagram
\begin{diagram}[75]
   B(G) & \Ear{\lambda_G} & T(\calO G) &&
   \Sar{\res^G_S} & & \saR{\res^G_S} &&
   B(S) & \Ear{\lambda_S} & T(\calO S) &&
\end{diagram}
and from (1) being the identity.
\end{proof}

For the proof of Theorem~C we need the following lemma which follows immediately from Lemma~\ref{lem Brauer construction 2}.

\begin{lemma}\label{lem beta and res}
Let $H\le G$. One has a commutative diagram 
\begin{diagram}[80]
  T(\calO G) & \movearrow(10,0){\Ear[40]{\beta_G}} & 
             \movevertex(60,-8){\mathop{\prod}\limits_{P\in\scrS_p(G)} R(K[N_G(P)/P])}  & & & &&
  \Sar{\res^G_H} & & \movearrow(60,0){\saR{\res^G_H}} & & & &&
  T(\calO H) & \movearrow(10,0){\Ear[40]{\beta_H}} & 
          \movevertex(60,-8){\prod\limits_{P\in\scrS_p(H)} R(K[N_H(P)/P])} & & & &&
\end{diagram}
where the right vertical map $\res^G_H$ is the ring homomorphism defined by sending $(\chi_P)_{P\in\scrS_p(G)}$ to 
$\bigl(\res^{N_G(P)/P}_{N_H(P)/P}(\chi_P)\bigr)_{P\in\scrS_p(H)}$.
\end{lemma}

\begin{proof} {\em of Theorem C.}\quad
First note that if we define $-^\circ\colon B(G)\to B(G)$ as the identity, then all four maps in (\ref{eqn splitting 2}) commute with $-^\circ$. Since they are also ring homomorphisms, they restrict to a sequence of group homomorphisms
\begin{equation}\label{eqn splitting 3}
   \lambda_S^{-1}\circ\res^G_S\circ\lambda_G\circ t_S^G\colon   
   B(\calF)^\times\to B(G)^\times \to \mathrm{O}(T(\calO G)) \to \mathrm{O}(T(\calO S)) = T(\calO S)^\times 
   \myiso B(S)^\times
\end{equation}
whose composition is the identity of $B(\calF)^\times$. It follows immediately that 
\begin{equation}\label{eqn decomp}
   \mathrm{O}(T(\calO G)) = (\lambda_G\circ t_S^G)(B(\calF))\times 
   \ker\bigr(\res^G_S\colon \mathrm{O}(T(\calO G))\to\mathrm{O}(T(\calO S))\bigr)
\end{equation}
and that the first factor is isomorphic to $B(\calF)^\times$, since $\lambda_G\circ t_S^G$ is injective. 

\smallskip
Next we show that $\beta_G\bigl(\ker\bigr(\res^G_S\colon \mathrm{O}(T(\calO G))\to\mathrm{O}(T(\calO S))\bigr)\bigr)$ is equal to the group
\begin{equation}\label{eqn Hom(-,K)}
    \Bigl(\prod_{P\in\scrS_p(G)} \Hom(N_G(P)/P,K^\times)\Bigr)'
 \end{equation}
consisting of all tuples $(\varphi_P)\in\bigl(\prod_{P\in\scrS_p(G)}\Hom(N_G(P)/P,K^\times)\bigr)^G$  satisfying
 \begin{equation}\label{eqn varphi condition}
    \varphi_P(xP) = \varphi_{P\langle x_p\rangle}(xP\langle x_P\rangle)\,,
 \end{equation}
for all $P\in\scrS_p(G)$ and $x\in N_G(P)$. First suppose that $u\in \mathrm{O}(T(\calO G))$ with $\res^G_S(u)=1=[\calO_S]$. Write $\beta_G(u)=(\epsilon_P\cdot\varphi_P)_{P\in\scrS_p(G)}$ as in Proposition~\ref{prop beta(u) decomposition}. Lemma~\ref{lem beta and res} applied to $H=S$ implies that $\epsilon_P=1$ for all $P\le S$. Since the tuple $(\epsilon_P)$ is a $G$-fixed point, we obtain $\epsilon_P=\epsilon_{(\lexp{g}{P})}$ for all $P\in\scrS_p(G)$ and $g\in G$. This implies that $\epsilon_P=1$ for all $P\in\scrS_p(G)$ and that $\beta_G(u)$ belongs to the group in (\ref{eqn Hom(-,K)}). Conversely, if a tuple $(\varphi_P)$ belongs to the group in (\ref{eqn Hom(-,K)}), then by Theorem~A there exists an element $\omega\in T(\calO G)$ with $\beta(\omega)=(\varphi_P)$. Since $\beta_G$ commutes with duals, we obtain $\beta_G(\omega^\circ)=(\varphi_P)^\circ=(\varphi_P^{-1})$ and therefore $\beta_G(\omega\cdot\omega^\circ)=\beta_G(\omega)\beta_G(\omega)^\circ = 1$. Since $\beta_G$ is injective, we obtain $\omega\in\mathrm{O}(T(\calO G))$.
 
 \smallskip
 Since $\beta_G$ is injective we have now proved that the second factor in (\ref{eqn decomp}) is isomorphic to the group in (\ref{eqn Hom(-,K)}). Finally, if $(\varphi_P)$ belongs to the group in (\ref{eqn Hom(-,K)}) then $\varphi_P$ vanishes on $p$-elements by the condition in (\ref{eqn varphi condition}) for all $P\in\scrS_p(G)$. Thus each $\varphi_P$ takes values in the roots of unity of $K$ of $p'$-order of $K$ which correspond isomorphically to those of $F^\times$.  This concludes the proof of Theorem~C. 
 \end{proof}

The next result is a further decomposition of the second factor occurring in the decomposition of $\mathrm{O}(T(\calO G))$ in Theorem~C.

\begin{proposition}\label{prop SO decomposition}
One has a decomposition
\begin{equation*}
   \Bigl(\prod_{P\in\scrS_p(G)} \Hom(N_G(P)/P,F^\times)\Bigr)'\cong \Hom(G,F^\times)\times
   \Bigl(\prod_{P\in\scrS_p(G)} \Hom(N_G(P)/PC_G(P),F^\times)\Bigr)'\,,
\end{equation*}
where the second factor denotes the set of all tuples 
\begin{equation*}
   (\varphi_P)\in \bigl(\prod_{P\in\scrS_p(G)} \Hom(N_G(P)/PC_G(P),F^\times)\bigr)^G
\end{equation*} 
satisfying 
\begin{equation*}
   \varphi_P(xPC_G(P))=\varphi_Q(xQC_G(Q))\,,
\end{equation*}
for all $P\in\scrS_p(G)$ and all $x\in N_G(P)$, where $Q:=P\langle x_p\rangle$.
\end{proposition}

\begin{proof}
Set
\begin{equation}\label{eqn varphi tuples}
   H:= \bigl(\prod_{P\in\scrS_p(G)} \Hom(N_G(P),F^\times)\bigr)^G
\end{equation}
and note that the component-wise inflations map the group $\bigl(\prod_{P\in\scrS_p(G)} \Hom(N_G(P)/P,F^\times)\bigr)'$ isomorphically onto the subgroup $H_1$ of $H$, consisting of those tuples $(\varphi_P)\in H$ that satisfy
\begin{equation}\label{eqn varphi condition 1}
   P\le \ker(\varphi_P)\quad \text{and}\quad \varphi_P(x)=\varphi_{P\langle x_p\rangle}(x)
\end{equation}
for all $P\in\scrS_p(G)$ and all $x\in N_G(P)$. Similarly, the component-wise inflations map the group 
$\bigl(\prod_{P\in\scrS_p(G)} \Hom(N_G(P)/PC_G(P),F^\times)\bigr)'$ isomorphically onto the group $H_2$ of tuples $(\varphi_P)\in H$ satisfying 
\begin{equation}\label{eqn varphi condition 2}
   PC_G(P)\le \ker(\varphi_P)\quad \text{and}\quad \varphi_P(x)=\varphi_{P\langle x_p\rangle}(x)
\end{equation}
for all $P\in\scrS_p(G)$ and all $x\in N_G(P)$.

\smallskip
Note that the homomorphisms $\iota\colon \Hom(G,F^\times)\to H_1$, $\varphi\mapsto (\varphi|_{N_G(P)})_{P\in\scrS_p(G)}$, and $\pi\colon H_1\to\Hom(G,F^\times)$, $(\varphi_P)\mapsto \varphi_1$, where $1\in\scrS_p(G)$ denotes the trivial subgroup of $G$, satisfy $\pi\circ\iota=\id$. Therefore, one has a decomposition
\begin{equation}\label{eqn final decomposition}
   \Bigl(\prod_{P\in\scrS_p(G)} \Hom(N_G(P)/P,F^\times)\Bigr)' \cong H_1 \cong \iota(\Hom(G,F^\times)) \times \ker (\pi)\,.
\end{equation}
Since $\iota$ is injective, the first factor in (\ref{eqn final decomposition}) is isomorphic to $\Hom(G,F^\times)$. Thus, it suffices to show that $\ker(\pi)=H_2$. Clearly, if $(\varphi_P)\in H_2$ then $\varphi_1=1$. Conversely, assume now that $(\varphi_P)\in H$ with $\varphi_1=1$. It suffices to show that $PC_G(P)\le \ker(\varphi_P)$ for all $P\in\scrS_p(G)$. Since every torsion element of $F^\times$ has $p'$-order, $\varphi_P$ is trivial on any $p$-element of $N_G(P)$. Thus, it suffices to show that $\varphi_P(c)=1$ for all $p'$-elements $c\in C_G(P)$. We show this by induction on $|P|$. The case $P=1$ is clear, since $\varphi_1=1$. Let $1<P\in\scrS_p(G)$ be arbitrary. Choose a normal subgroup $Q$ of $P$ of index $p$ and an element $x\in P$ with $P=Q\langle x\rangle$. Then, since $\varphi_P$ and $\varphi_Q$ vanish on $p$-elements and by the second property in (\ref{eqn varphi condition 1}) applied to $Q$ and $xc$, we obtain $\varphi_P(c)=\varphi_P(xc)=\varphi_Q(xc)=\varphi_Q(c)$ which equals $1$ by induction.
\end{proof}


\section{A criterion \`a la Yoshida}\label{sec Yoshida criterion}

In this section we prove a result analogous to Yoshida's characterization of the mark values of units of the Burnside ring, see \cite[Proposition~6.5]{Yoshida1990}. To set the stage, recall from Remark~\ref{rem orthogonal units}(b) that one has an injective ring homomorphism
\begin{equation*}
   \sigma_G\colon T(\calO G)\to \Bigl(\prod_{(P,sP)\in\scrT_p(G)} K\Bigr)^G\,,\quad 
   \omega\mapsto \bigl(  \kappa_{N_G(P)/P}(\omega(P))(sP)  \bigr)_{(P,sP)\in\scrT_p(G)}\,,
\end{equation*}
where we view the element $\kappa_{N_G(P)/P}(\omega(P))\in R(K[N_G(P)/P])$ as virtual character that can be evaluated at $sP$, and where the exponent $G$ stands for taking fixed points under the natural conjugation action $\lexp{g}{\bigl(a_{(P,sP)}\bigr)_{(P,sP)\in\scrT_p(G)}}:=\bigl(\lexp{g}{a_{\lexp{g^{-1}}{(P,sP)}}}\bigr)_{(P,sP)\in\scrT_p(G)}$. Here we use the action of $G$ on $\scrT_p(G)$ given by $\lexp{g}{(P,sP)}=(\lexp{g}{P},\lexp{g}{(sP)})$.
Since every orthogonal unit of $T(\calO G)$ has finite order (see Remark~\ref{rem orthogonal units}(b) or Theorem~C), the ring homomorphism $\sigma_G$ restricts to an injective group homomorphism
\begin{equation*}
   \sigma_G\colon \mathrm{O}(T(\calO G))\to \Bigl(\prod_{(P,sP)\in\scrT_p(G)} \langle\zeta\rangle\Bigr)^G\,,
\end{equation*}
where $\zeta\in K^\times$ is a root of unity that generates the group of roots of unity in $K$ of order dividing $2\exp(G)_{p'}$ (see again Theorem~C).

\begin{theorem}\label{thm a la Yoshida}
Let $\zeta\in K^\times$ be as above and let $(\alpha_{(P,sP)})\in\bigl(\prod_{(P,sP)\in\scrT_p(G)} \langle\zeta\rangle\bigr)^G$.
The following are equivalent:

\smallskip
{\rm (i)} The tuple $(\alpha_{(P,sP)})$ belongs to $\sigma_G(\mathrm{O}(T(\calO G)))$.

\smallskip
{\rm (ii)} For every $P\in\scrS_p(G)$, the function $\psi_P\colon N_G(P)/P\to\langle \zeta \rangle$, 
$xP\mapsto \alpha_{(Q, xQ)}\cdot\alpha_{(P,1P)}$, where $Q=P\langle x_p \rangle$, is a group homomorphism.
\end{theorem}

\begin{proof}
Note that the functions $\psi_P$ in (ii) are well-defined: If $xP=yP$ for $x,y\in N_G(P)$, then $P\langle x_p\rangle=P\langle y_p \rangle$ and $xP\langle x_p\rangle = yP\langle y_p\rangle$.

Suppose first that $u\in\mathrm{O}(T(\calO G))$ and $(\alpha_{(P,sP)}) =\sigma_G(u)$. By Proposition~\ref{prop beta(u) decomposition} we have  $\beta_G(u)=(\epsilon_P\cdot\varphi_P)$ with $\epsilon_P\in\{\pm 1\}$ and $\varphi_P\in\Hom(N_G(P)/P, K^\times)$. Thus, $\alpha_{(P,sP)} = \epsilon_P\cdot \varphi_P(sP)$ for all $(P,sP)\in\scrT_p(G)$. As a consequence, for all $P\in\scrS_p(G)$ and $xP\in N_G(P)/P$, we have $\psi_P(xP)= \epsilon_Q\cdot\varphi_Q(xQ)\cdot\epsilon_P$, where $Q=P\langle x_p\rangle$. 
However, by the coherence property of the components of $\beta_G(u)$, we have $\epsilon_Q\cdot\varphi_Q(xQ)=\epsilon_P\cdot\varphi_P(xP)$ so that $\psi_P(xP)=\varphi_P(xP)$ and $\psi_P=\varphi_P$ is a homomorphism.

\smallskip
Conversely, suppose that $\psi_P$ is a homomorphism for every $P\in\scrS_p(G)$. This implies $1=\psi_P(1P) = \alpha_{(P,1P)}\cdot\alpha_{(P,1P)}$ so that $\epsilon_P:=\alpha_{(P,1P)}\in\{\pm1\}$. For $P\in\scrS_p(G)$, we define $\chi_P\colon N_G(P)/P\to K$, $xP\mapsto \epsilon_P\cdot\psi_P(xP)$. Since $\psi_P\in\Hom(N_G(P)/P,K^\times)$, $\chi_P=\epsilon_P\cdot\psi_P\in R(K[N_G(P)/P])$. Since the tuple $(\alpha_{(P,sP)})$ is $G$-fixed, the tuple $(\chi_P)$ is $G$-fixed. We check that the tuple $(\chi_P)$ satisfies the coherence condition. Let $P\in\scrS_p(G)$, $x\in N_G(P)$ and set $Q:=P\langle x_p\rangle$. Then $\chi_P(xP)=\epsilon_P\cdot\psi_P(xP)=\epsilon_P^2\cdot\alpha_{(Q,xQ)} = \alpha_{(Q,xQ)}$ and also $\chi_Q(xQ)=\epsilon_Q\cdot\psi_Q(xQ) = \epsilon_Q^2\cdot\alpha_{(Q,xQ)}=\alpha_{(Q,xQ)}$. By Theorem~A, there exists $u\in T(\calO G)$ with $\beta_G(u)=(\chi_P)$, which implies that the $(P,sP)$-component of $\sigma_G(u)$ is $\chi_P(sP)=\alpha_{(P,sP)}$, so that $\sigma_G(u)=(\alpha_{(P,sP)})$. Moreover, $\beta_G(u\cdot u^\circ)=\beta_G(u)\cdot\beta_G(u)^\circ = 1$, since $\beta_G(u)^\circ=(\chi_P)^\circ = (\epsilon_P\cdot\psi_P^{-1})$. The injectivity of $\beta_G$ implies $u\cdot u^\circ=1$ so that $u\in \mathrm{O}(T(\calO G))$.
\end{proof}


\section{Examples}\label{sec examples}

\begin{proposition}
Suppose that every $p$-singular element of $G$ is a $p$-element. Let $S$ be a Sylow $p$-subgroup of $G$ and let $\calF:=\calF_S(G)$ be the associated fusion system on $S$. Moreover, let $\scrStilde_p(G)$ be a set of representatives of the $G$-conjugacy classes of $p$-subgroups of $G$. Then
\begin{equation*}
   \mathrm{O}(T(\calO G)) \cong B(\calF)^\times\times\prod_{P\in\scrStilde_p(G)} \Hom(N_G(P)/P,F^\times)\,.
\end{equation*}
\end{proposition}

\begin{proof}
We use the decomposition of $\mathrm{O}(T(\calO G))$ from Theorem~C. Since every $p$-singular element of $G$ is a $p$-element, the coherence condition in (\ref{eqn varphi coherence condition}) is satisfied for every tuple $(\varphi_P)\in\bigl(\prod_{P\in\scrS_p(G)} \Hom(N_G(P)/P,F^\times)\bigr)^G$. Thus,
\begin{equation*}
  \Bigl(\prod_{P\in\scrS_p(G)} \Hom(N_G(P)/P,F^\times)\Bigr)' =  \Bigl(\prod_{P\in\scrS_p(G)} \Hom(N_G(P)/P,F^\times)\Bigr)^G\,.
\end{equation*}
Moreover, the natural projection map
\begin{equation*}
   \Bigl(\prod_{P\in\scrS_p(G)} \Hom(N_G(P)/P,F^\times)\Bigr)^G \to \prod_{P\in\scrStilde_p(G)} \Hom(N_G(P)/P,F^\times)
\end{equation*}
is an isomorphism.
\end{proof}

Note that the hypothesis of the above proposition is satisfied if $G$ is a Frobenius group of the form $S\rtimes E$, where the Frobenius kernel $S$ is a $p$-group and $E\le\Aut(S)$ is a $p'$-group. In particular, if $S$ is cyclic, we have the following corollary.

\begin{corollary}
Suppose that $G=S\rtimes E$, where $S$ is a cyclic $p$-group of order $p^n$ and $E\le \Aut(S)$ is a $p'$-group. Then
\begin{equation*}
   \mathrm{O}(T(\calO G)) \cong B(S)^\times\times\prod_{P\in\scrStilde_p(G)} \Hom(E,F^\times)^{n+1}\,.
\end{equation*}
\end{corollary}

\begin{proof}
This follows immediately from the previous proposition noting that for every $P\in\scrS_p(G)$ one has $N_G(P)=G$ and $\Hom(G/P,F^\times)\cong \Hom(E,F^\times)$. Moreover, $B(\calF)=B(S)$, since $S$ is cyclic.
\end{proof}

Recall that $G$ is called {\em $p$-nilpotent} if there exists a normal $p'$-subgroup $N$ of $G$ such that $G/N$ is a $p$-group.

\begin{proposition}
Suppose that $G$ is a $p$-nilpotent group with Sylow $p$-subgroup $S$ and let $\calF:=\calF_S(G)$ be the associated fusion system on $S$. Then
\begin{equation*}
   \mathrm{O}(T(\calO G)) \cong B(S)^\times \times \Hom(G,F^\times)\,.
\end{equation*}
\end{proposition}

\begin{proof}
We use again Theorem~C. Since $G$ is $p$-nilpotent, one has $\calF_S(G)=\calF_S(S)$ by a theorem of Frobenius. Therefore, two subgroups $P$ and $Q$ of $S$ are $G$-conjugate if and only if they are $S$-conjugate. It follows that $B(\calF)=B(S)$. Moreover, since $G$ is $p$-nilpotent, the group $N_G(P)/C_G(P)$ is a $p$-group for every $P\in\scrS_p(G)$. The result now follows from Proposition~\ref{prop SO decomposition} and noting that $\Hom(N_G(P)/PC_G(P),F^\times)$ is trivial for all $P\in\scrS_p(G)$.
\end{proof}


\end{document}